\documentclass[10pt]{article}

\usepackage{amscd,amsmath, amssymb, fancyhdr, mathbbol}

\numberwithin{equation}{section}

\def\eqref#1{(\ref{#1})}
\newcommand{\goth}{\mathfrak}

\newcommand{\arrow}{{\:\longrightarrow\:}}
\newcommand{\Z}{{\Bbb Z}}
\def\C{{\Bbb C}}

\newcommand{\Q}{{\Bbb Q}}

\def\1{\sqrt{-1}\:}

\newcommand{\cntrct}                
{\hspace{2pt}\raisebox{1pt}{\text{$\lrcorner$}}\hspace{2pt}}

\renewcommand{\tilde}{\widetilde}
\renewcommand{\bar}{\overline}
\renewcommand{\phi}{\varphi}
\renewcommand{\epsilon}{\varepsilon}
\renewcommand{\geq}{\geqslant}
\renewcommand{\leq}{\leqslant}

\newcommand{\Teich}{\operatorname{Teich}}

\newcommand{\Aut}{\operatorname{Aut}}

\newcommand{\Diff}{\operatorname{Diff}}

\newcommand{\Comp}{\operatorname{Comp}}

\newcounter{Mycounter}[section]
\newcounter{lemma}[section]
\renewcommand{\thelemma}{{Lemma \thesection.\arabic{lemma}}}
\newcommand{\lemma}{%
    \setcounter{lemma}{\value{Mycounter}}
    \refstepcounter{lemma}
    \stepcounter{Mycounter}
    {\noindent \bf \thelemma:\ }}

\newcounter{claim}[section]
\newcounter{sublemma}[section]
\newcounter{corollary}[section]
\renewcommand{\thecorollary}{{Corollary \thesection.\arabic{corollary}}}
\newcommand{\corollary}{%
    \setcounter{corollary}{\value{Mycounter}}
    \refstepcounter{corollary}
    \stepcounter{Mycounter}
    {\noindent \bf \thecorollary:\ }}

\newcounter{theorem}[section]
\renewcommand{\thetheorem}{{Theorem \thesection.\arabic{theorem}}}
\newcommand{\theorem}{%
    \setcounter{theorem}{\value{Mycounter}}
    \refstepcounter{theorem}
    \stepcounter{Mycounter}
    {\noindent \bf \thetheorem:\ }}

\newcounter{conjecture}[section]
\renewcommand{\theconjecture}{{Conjecture \thesection.\arabic{conjecture}}}
\newcommand{\conjecture}{%
    \setcounter{conjecture}{\value{Mycounter}}
    \refstepcounter{conjecture}
    \stepcounter{Mycounter}
    {\noindent \bf \theconjecture:\ }}

\newcounter{proposition}[section]
\newcounter{definition}[section]
\renewcommand{\thedefinition}
      {{Definition~\thesection.\arabic{definition}}}
\newcommand{\definition}{%
    \setcounter{definition}{\value{Mycounter}}
    \refstepcounter{definition}
    \stepcounter{Mycounter}
    {\noindent \bf \thedefinition:\ }}

\newcounter{example}[section]
\newcounter{remark}[section]
\renewcommand{\theremark}{{Remark \thesection.\arabic{remark}}}
\newcommand{\remark}{%
    \setcounter{remark}{\value{Mycounter}}
    \refstepcounter{remark}
    \stepcounter{Mycounter}
    {\noindent \bf \theremark:\ }}

\newcounter{problem}[section]
\newcounter{question}[section]
\makeatletter

\@addtoreset{equation}{section} \@addtoreset{footnote}{section}
\makeatother

\def\blacksquare{\hbox{\vrule width 5pt height 5pt depth 0pt}}
\def\endproof{\blacksquare}

\begin{document}
\begin{center}
{\LARGE\bf
Survey of finiteness results for hyperk\"ahler manifolds
\\[4mm]
}

Ljudmila Kamenova
\footnote{This work was partially supported by grant 346300 for IMPAN from 
the Simons Foundation and the matching 2015-2019 Polish MNiSW fund}

\end{center}

{\small \hspace{0.1\linewidth}
\begin{minipage}[t]{0.8\linewidth}
{\bf Abstract} \\
This paper is a survey of finiteness results in hyperk\"ahler geometry. 
We review some classical theorems by Sullivan, Koll\'ar-Matsusaka, Huybrechts, 
as well as theorems in the recent literature by Charles, Sawon, 
and joint results of the author with Verbitsky. We also strenghten a 
finiteness theorem of the author. 
These are extended notes of the author's talk during the closing conference 
of the Simons Semester in the Banach Center in B\c{e}dlewo, Poland. 

\end{minipage}
}

{\scriptsize
\tableofcontents
}


\section{Introduction}


In each given complex dimension $2n$ there are only a few known examples 
of compact hyperk\"ahler manifolds up to deformation, namely, the Hilbert 
scheme of $n$ points on a K3 surface $S$, the generalized Kummer variety  
$K^n(A)$ and the exceptional examples $O_6$ (for $n=3$) and $O_{10}$ 
(for $n=5$) given by O'Grady. A natural question to ask is if there are 
only finitely many compact hyperk\"ahler manifolds in any given dimension,  
up to deformation. This paper surveys some of the known results in this 
direction. 

\hfill

Based on Koll\'ar-Matsusaka's finiteness theorem, Daniel Huybrechts proved 
the following theorem in \cite{_Huybrechts:finiteness_}. 

\hfill

\theorem (Huybrechts, \cite{_Huybrechts:finiteness_}) \label{H1}
If the second integral cohomology group $H^2(\Z)$ and the homogeneous polynomial
of degree $2n-2$ on $H^2(\Z)$ defined by the first Pontryagin class are 
given, then there exist at most finitely many diffeomorphism types of 
compact hyperk\"ahler manifolds of complex dimension $2n$ realizing this 
structure.

\hfill

Huybrechts \cite{_Huybrechts:finiteness_} also showed that fixing the  
Beauville-Bogomolov-Fujiki form (i.e., giving the abelian group $H^2(\Z)$ 
a ring structure) is equivalent to fixing the the first Pontryagin class. 

\hfill

If the diffeomorphic structure on a manifold $M$ is given, then there are also 
finitely many deformation types of hyperk\"ahler metrics on $M$. 

\hfill

\theorem (Huybrechts, \cite{_Huybrechts:finiteness_}) \label{H2}
Let $M$ be a fixed compact manifold. Then there exist at most finitely
many different deformation types of irreducible holomorphic symplectic 
complex structures on $M$.

\hfill 

Using \ref{H2}, the author and Misha Verbitsky established the following 
finiteness results in \cite{_KV:fibrations_} about hyperk\"ahler fibrations. 

\hfill

\theorem (Kamenova-Verbitsky, \cite{_KV:fibrations_})
Let $M$ be a fixed compact manifold of complex dimension $2n$ and 
$b_2(M) \geq 7$. Then there are 
only finitely many deformation types of hyperk\"ahler 
Lagrangian fibrations $M \arrow \C P^n$. 

\hfill

Fran\c{c}ois Charles has the following boundedness result for families of 
hyperk\"ahler varieties up to deformation. This sharpens Koll\'ar-Matsusaka's 
finiteness result for hyperk\"ahler manifolds, because he replaces the 
assumption that $L$ is ample with the weaker assumption $q(L)>0$, but one must 
still bound the top self-intersection of $L$. 

\hfill 

\begin{theorem}(Charles, \cite{_Charles_})
Let $n$ and $r$ be two positive integers. Then there exists a scheme $S$ of 
finite type over $\C$, and a projective morphism ${\cal M} \arrow S$ such 
that if $M$ is a complex hyperk\"ahler variety of dimension $2n$ and $L$ is a 
line bundle on $M$ with $L^{2n} = r$ and $q(L) > 0$, where $q$ is the
Beauville-Bogomolov form, then there exists a complex point $s$ of $S$ 
such that ${\cal M}_s$ is birational to $M$. 
\end{theorem}

\hfill

For a hyperk\"ahler manifold $M$, the Fujiki constant and the 
discriminant of the Beauville-Bogomolov-Fujiki form $q$ are topological 
invariants. If we fix them instead of fixing the whole intersection 
form $q$, we can ask for finiteness given much less data than in 
\ref{H1}. Here we establish the following result, which generalizes our 
theorem in \cite{_Kamenova_} previously stated for Lagrangian fibrations. 

\hfill 

\theorem
There are at most finitely many deformation classes of hyperk\"ahler manifolds 
with a fixed Fujiki constant $c$ and a given 
discriminant of the Beauville-Bogomolov-Fujiki lattice $(\Lambda,q)$. 

\hfill

In our original proof we used a primitive vector $v$ with $q(v)=0$ coming 
from the line bundle $L$ associated to the Lagrangian fibration. However, 
we do not necessarily need $v$ to come from a Lagrangian fibration. 

\hfill

In \cite{_Sawon_finit_}, Sawon proved a finiteness result for Lagrangian 
fibrations with several assumptions on the fibration. We'll formulate Sawon's 
theorem in section \ref{C_S}. In \cite{_Kamenova_} we gave the following 
generalization of his result. 

\hfill

\theorem (Kamenova, \cite{_Kamenova_})
Consider a Lagrangian fibration $\pi:M \arrow\C P^n$ such that 
there is a line bundle $P$ on $M$ with $q(P)>0$ and with a given $P$-degree 
$d$ on the general fiber $F$ of $\pi$, i.e., $P^n \cdot F = d$. 
Then there are at most finitely many deformation classes of hyperk\"ahler 
manifolds $M$ as above, i.e., they form a bounded family.


\section{Basic results in hyperk\"ahler geometry}


\definition  A {\bf hyperk\"ahler manifold}
is a compact simply connected K\"ahler holomorphic symplectic manifold. 
A hyperk\"ahler manifold $M$ is called {\bf irreducible} if $H^{2,0}(M)=\C$.

\hfill

According to Bogomolov's decomposition theorem,
\cite{_Bogomolov:decompo_}, 
any hyperk\"ahler manifold admits a finite covering 
which is a product of finitely many irreducible hyperk\"ahler manifolds. 
From now on we shall assume that all hyperk\"ahler manifolds are irreducible. 

\hfill

\remark 
In the compact case the following two notions are equivalent: 
a holomorphic symplectic K\"ahler manifold and 
a manifold with a {\em hyperk\"ahler structure}, 
that is, a triple of complex structures 
satisfying the quaternionic relations and 
parallel with respect to the Levi-Civita 
connection. This equivalence follows from Yau's solution 
of Calabi's conjecture (\cite{_Besse:Einst_Manifo_}). 
Throughout this paper we assume compactness and we use the 
complex algebraic point of view. 

\hfill

\definition \label{monodromy-def}
Let $M$ be a compact complex manifold and $\Diff^0(M)$ the connected 
component of the identity of its diffeomorphism group. 
Denote by $\Comp$ the space of complex structures on $M$, equipped with 
a structure of Fr\'echet manifold. The {\em Teichm\"uller space} of $M$ 
is the quotient  $\Teich:=\Comp/\Diff^0(M)$. 
For a hyperk\"ahler manifold $M$, the Teichm\"uller space is 
finite-dimensional (\cite{_Catanese:moduli_}). 
Let $\Diff^+(M)$ be the group of orientable diffeomorphisms of 
a complex manifold $M$. The {\em mapping class group} 
$$\Gamma := \Diff^+(M)/\Diff^0(M)$$ acts naturally on $\Teich$. 
For $I\in \Teich$, let $\Gamma_I$ be the subgroup of $\Gamma$ which fixes 
the connected component $\Teich_I$ of the complex structure $I$. 
The {\em monodromy group} is the image of $\Gamma_I$ in $\Aut H^2(M, \Z)$. 

\hfill

On $H^2(M, \Z)$ there is a natural primitive integral quadratic form, 
called the {\bf  Beauville-Bogomolov-Fujiki form}, or {\bf BBF form} 
for shortness. The easiest way 
to define it is via the Fujiki relation below. For the classical 
definition we refer the reader to \cite{_Beauville_} and 
\cite{_Huybrechts:lec_}. 

\hfill

\theorem
(Fujiki, \cite{_Fujiki:HK_}) \label{Fujiki_formula}
Let $\eta\in H^2(M, \Z)$ and $\dim M=2n$, where $M$ is a 
hyperk\"ahler manifold. Then $\int_M \eta^{2n}= c \cdot q(\eta,\eta)^n$,
for a primitive integral quadratic form $q$ on $H^2(M, \Z)$, where $c>0$ is a 
constant depending on the topological type of $M$. The constant $c$ in Fujiki's 
formula is called the {\bf Fujiki constant}. 
\endproof

\hfill



\remark
The form $q$ has signature $(3,b_2-3)$. 
It is negative definite on primitive forms and positive
definite on the space $\langle \Omega, \bar \Omega, \omega\rangle$, 
where $\Omega$ is the holomorphic symplectic form 
and $\omega$ is a K\"ahler form 
(see \cite{_Verbitsky:cohomo_}, Theorem 6.1 
and \cite{_Huybrechts:lec_}, Corollary 23.9).

\hfill

\definition
Let $\eta\in H^{1,1}(M)$ be a real (1,1)-class on
a hyperk\"ahler manifold $M$. We say that $\eta$
is {\bf parabolic} if $q(\eta,\eta)=0$.
A line bundle $L$ is called {\bf parabolic} if the class $c_1(L)$
is parabolic.

\hfill

\remark \label{_P_L_identity_}
If $L$ is a parabolic class and $P \in H^2(M)$ is any class, then 
after we substitute $\eta = P+tL$ into Fujiki's formula in 
\ref{Fujiki_formula}, and compare the coefficients of $t^n$ on both sides, 
we obtain ${2n \choose n} P^n L^n = c 2^n q(P,L)^n$. 

\hfill

Notice that hyperk\"ahler manifolds have a very restricted fibration 
structure. 

\hfill

\theorem
(Matsushita, \cite{_Matsushita:fibred_}).
Let $\pi:\; M \rightarrow B$ be a surjective holomorphic map with 
connected fibers 
from a hyperk\"ahler manifold $M$ to a base $B$, with $0<\dim B < \dim M$.
Then $\dim B = \frac{1}{2} \dim M$, and the fibers of $\pi$ are 
holomorphic Lagrangian (i.e., the symplectic
form vanishes when restricted to the fibers).

\hfill

Such a map is called
{\bf a holomorphic Lagrangian fibration}.

\hfill

\remark 
D. Matsushita (\cite{_Matsushita:CP^n_}) proved that if the base of 
$\pi$ is smooth and $M$ is projective, $B$ has the same rational cohomology 
as $\C P^n$. Later J.-M. Hwang (\cite{_Hwang:base_}) proved that 
under the same assumptions $B\cong \C P^n$. 

\hfill

\definition A line bundle $L$ is called 
{\bf semiample} if  $L^N$ is generated 
by its holomorphic sections which have  
no common zeros. 

\hfill



\remark \label{_main_rem}
From semiampleness 
it trivially follows that $L$ is nef, however a nef bundle 
is not necessarily semiample. 
Let $\pi:\; M \rightarrow B$ 
be a holomorphic Lagrangian fibration and let $\omega_B$ be 
a K\"ahler class on $B$. Then $\eta:=\pi^*\omega_B$ is 
semiample and parabolic. By Matsushita's theorem, 
the converse is also true:
if $L$ is semiample and parabolic, $L$ induces a Lagrangian
fibration. 

\hfill

\conjecture\label{_SYZ_conj_Conjecture_}
({\bf Hyperk\"ahler SYZ conjecture})
Let $L$ be a parabolic nef line bundle
on a hyperk\"ahler manifold. Then
$L$ is semiample.

\hfill

\remark
The SYZ conjecture can be seen as
a hyperk\"ahler version of the ``abundance conjecture''
(see e.g. \cite{_Demailly_Peternell_Schneider:ps-eff_}, 
2.7.2). It was stated by many authors such as Tyurin, Bogomolov, 
Hassett-Tschinkel, Huybrechts, Sawon, Verbitsky, etc. For more details on 
the SYZ conjecture one might look into \cite{_Sawon_ab_} and 
\cite{_Verbitsky:SYZ_}. 

\hfill

\remark \label{SYZ-Lag}
As a corollary to the SYZ conjecture we see that any hyperk\"ahler 
manifold $M$ with $b_2(M) \geq 5$ can be deformed to one that admits a 
Lagrangian fibration. This is one of the reasons to study Lagrangian 
fibrations on hyperk\"ahler manifolds. Indeed, if $b_2(M) \geq 5$, 
by Meyer's theorem there is an isotropic vector $v$. One can deform 
$M$ to a hyperk\"ahler manifold $M'$ whose Picard group is generated by 
a primitive line bundle $L$ corresponding to $v$, i.e., $q(c_1(L))=0$. 
The positive cone of $M'$ coincides with the K\"ahler cone, and therefore, 
either $L$ or $L^*$ is nef. By the SYZ conjecture, $M'$ admits a 
Lagrangian fibration. This is the argument of Proposition 4.3 in 
Sawon's paper \cite{_Sawon_ab_}. 

\hfill

Relatively recently Matsushita conjectured that the fibration structure of a 
hyperk\"ahler manifold is even more restricted than the structure 
given by Matsushita's theorem. The following conjecture was introduced 
to the author in private communictions 
with J. Sawon, D. Matsushita and J.-M. Hwang in 2013/14. 

\hfill

\conjecture\label{_Matsushita_Conjecture_}
({\bf Matsushita's conjecture}) 
Every holomorphic Lagrangian fibration $\pi:M\rightarrow\C P^n$ is 
either locally isotrivial or the fibers vary maximally in the 
moduli space of Abelian varieties ${\cal A}_n$. 

\hfill

B. van Geemen and C. Voisin recently proved the following weeker version of 
Matsushita's conjecture. Here $\rho$ is the rank of the Picard group. 

\hfill

\begin{theorem} (B. van Geemen, C. Voisin, \cite{_vG_V:matsushita_}) 
\label{gv-m}
Let $X$ be a projective hyperk\"ahler manifold of dimension $2n$ admitting a 
Lagrangian fibration $f:X\rightarrow B$. 
Assume $b_{2, tr}(X) = b_2(X) - \rho (X) \geq 5$. Then a very general deformation 
$(X',f',B')$ of the triple $(X,f,B)$ satisfies Matsushita's conjecture. 
\end{theorem}


\section{Some classical results}


There is a general classical result of Dennis Sullivan about finiteness of 
simply connected K\"ahler manifolds. It follows from the $\Q$-formality of 
K\"ahler manifolds combined with the description of the rational homotopy type 
of K\"ahler manifolds in terms of a differential graded algebra. 
This algebra encodes the cohomology with its ring structure as a 
primary invariant. 

\hfill

\begin{theorem}(Sullivan, \cite{_Sullivan_}) 
The diffeomorphism type of a simply connected K\"ahler manifold (of 
complex dimension greater than $2$) is finitely determined by its 
integral cohomology ring and its Pontryagin classes. 
\end{theorem}

\hfill

Daniel Huybrechts noticed that for a hyperk\"ahler manifold $M$ much less of 
the topology needs to be fixed in order to determine the diffeomorphism 
type up to finite ambiguity. The first Pontryagin class $p_1(M) \in 
H^4(M, \Z)$ gives rise to a homogeneous polynomial 
$\tilde p_1: H^2(M, \Z) \arrow \Z$ of degree $2n-2$.

\hfill

\theorem (Huybrechts, \cite{_Huybrechts:finiteness_}) 
If the second integral cohomology group $H^2(\Z)$ and the homogeneous polynomial
of degree $2n-2$ on $H^2(\Z)$ defined by the first Pontryagin class are 
given, then there exist at most finitely many diffeomorphism types of 
compact hyperk\"ahler manifolds of complex dimension $2n$ realizing this 
structure.

\hfill

Here Huybrechts fixes $H^2(\Z)$ as an abelian group. The first Pontryagin 
class, or equivalently, the Beauville-Bogomolov-Fujiki form, give the 
ring structure of the second cohomology $H^2(\Z)$. 
The proof of this theorem uses Koll\'ar-Matsusaka's finiteness result: 

\hfill

\begin{theorem} (Koll\'ar-Matsusaka, \cite{_KM:finiteness_}) \label{KM-f}
There are finitely many deformation types of projective 
manifolds $M$ of dimension $d$ that admit an ample line bundle $L$ 
with fixed intersection numbers $L^d$ and $K_M \cdot L^{d-1} \in \Z$.
\end{theorem}

\hfill

In particular, for Calabi-Yau manifolds $M$ of dimension $d$ it is 
enough to fix the top intersection $L^d$. In order to remove the 
fixed polarization, Huybrechts shows that any hypek\"ahler 
manifold with fixed $H^2(\Z)$ and $\tilde p_1$ can be deformed to one 
that admits a polarization $L$ with bounded $L^d$, using 
Hitchin-Sawon's formulas in \cite{_H_S_}. Here is a version of 
Huybrechts's finiteness 
theorem where the BBF form is given in place of $\tilde p_1$. 

\hfill 

\theorem (Huybrechts, \cite{_Huybrechts:finiteness_}) 
There are only finitely many deformation types of hyperk\"ahler 
manifolds $M$ of fixed dimension such that the lattice 
$(H^2(M, \Z), q_M)$ is isomorphic to a given one. 

\hfill

Once the diffeomorphism type of $M$ is fixed, Huybrechts also 
shows that there are finitely many deformation types of 
hyperk\"ahler metrics $g$ on $M$. 

\hfill

\theorem (Huybrechts, \cite{_Huybrechts:finiteness_}) \label{H-f-3}
Let $M$ be a fixed compact manifolds. Then there exist at 
most finitely many deformation types of hyperk\"ahler 
structures on $M$. 

\hfill

The idea of the proof is the following. 
Let $\Lambda = (H^2(M, \Z), q)$ and ${\goth M}_\Lambda$ be the 
coarse moduli space of marked hyperk\"ahler manifolds $(M, \phi)$. 
Fix a primitive positive element $v \in H^2(M, \Z)$, i.e., $q(v)>0$. 
An isomorphism, or a marking, $\phi : (H^2(M, \Z), q) \arrow \Lambda$ 
determines a point $(M, \phi) \in {\goth M}_\Lambda$. By Huybrechts's 
proectivity criterion \cite{_Huybrechts:basic_}, since $q(v)>0$,  
$M$ is deformation equivalent to a projective manifold $M'$ with an 
ample line bundle $L$ corresponding to $v$. By Fujiki's formula in 
\ref{Fujiki_formula}, $L^{2n} = c \cdot q(L)^n = c \cdot q(v)^n$. 
Fujiki's constant $c$ depends only on the topological type of $M$, 
which is fixed. Since $v$ is also fixed, $q(v)$ is determined, and 
therefore the top intersection $L^{2n}$ is known and one applies 
Koll\'ar-Matsusaka's \ref{KM-f}.


\section{Some recent results} \label{C_S}


Together with Misha Verbitsky in \cite{_KV:fibrations_} 
we were interested in studying finiteness 
questions about hyperk\"ahler manifolds that admit Lagrangian fibrations. 
As we pointed out in \ref{SYZ-Lag}, if the SYZ conjecture holds, then 
any hyperk\"ahler manifold (with $b_2 \geq 5$) can be deformed to one that 
admits a Lagrangian fibration, and therefore it would be enough to study 
Lagrangian fibrations for questions concerning deformation classes 
of hyperk\"ahler varieties. By Huybrechts's \ref{H-f-3}, for a fixed 
compact manifold there are at most finitely many deformation types of 
hyperk\"ahler structures on it. We proved that for each hyperk\"ahler 
structure there are only finitely many ways in which it fibers over a 
smooth base. 

\hfill

\theorem (Kamenova-Verbitsky, \cite{_KV:fibrations_})  \label{KVfin}
Let $M$ be a fixed compact manifold of complex dimension $2n$ with 
$b_2(M) \geq 7$. Then there are 
only finitely many deformation types of hyperk\"ahler 
Lagrangian fibrations $p: M \arrow \C P^n$. 

\hfill

Let $\Teich$ be the Teichm\"uller space of $M$ and $\Gamma_I < \Gamma$ 
as in \ref{monodromy-def}. If $M$ admits a Lagrangian fibration, there 
is a natural nef line bundle $L=p^* {\cal O} (1)$ associated to the fibration.  
Consider the set 
$$\Teich_L = \{ I \in \Teich | L \in H^{1,1}((M, I), \Z) \},$$ where $L$ 
remains of type $(1,1)$ on deformations of the complex structure. 
In \cite{_KV:fibrations_} 
we proved that there are finitely many orbits of the action of $\Gamma_I$ 
on $\Teich_L$ using lattice theory methods and Nikulin-style 
techniques applied to the BBF form. If we take a general deformation of 
a Lagrangian fibration, the Picard group would be 1-dimensional and generated 
by the nef line bundle $L$. For each such pair $(M,L)$ we prove that there 
is a unique deformation type of a fibration structure. 
Since there are finitely many orbits of $\Gamma_I$ acting on $\Teich_L$ 
we conclude finiteness of the deformation types of Lagrangian fibrations.

\hfill

Recently, Fran\c{c}ois Charles sharpened Koll\'ar-Matsusaka's theorem 
for hyperk\"ahler varieties by replacing the ampleness assumption with a 
weaker one. Our results, \ref{_finiteness_Theorem_1_} and 
\ref{_finiteness_Theorem_2_}, rely on his theorem. 

\hfill

\begin{theorem}(Charles, \cite{_Charles_}) \label{charles}
Let $n$ and $r$ be two positive integers. Then there exists a scheme $S$ of 
finite type over $\C$, and a projective morphism ${\cal M} \arrow S$ such 
that if $M$ is a complex hyperk\"ahler variety of dimension $2n$ and $L$ is a 
line bundle on $M$ with $c_1(L)^{2n} = r$ and $q(L) > 0$, where $q$ is the
Beauville-Bogomolov form, then there exists a complex point $s$ of $S$ 
such that ${\cal M}_s$ is birational to $M$. 
\end{theorem}

\hfill


Consider a lattice $\Lambda$, i.e., a free $\Z-$module of finite rank 
together with a non-degenerate symmetric bilinear from $q$ with 
values in $\Z$. If $\{ e_i \}$ is a basis of $\Lambda$, the 
{\it discriminant} of $\Lambda$ is $\text{discr}(\Lambda) = 
\text{det}(q(e_i,e_j))$. Let us recall the following lemma 
from \cite{_Kamenova_}. 

\hfill

\lemma\label{_bounded_lemma_}
Let $(\Lambda, q)$ be an indefinite lattice and $v \in \Lambda$ be an 
isotropic non-zero vector. Then there exists a positive vector 
$w \in \Lambda$ such that $0 < q(w,v) \leq |\text{discr}(\Lambda)|$ 
and $0<q(w,w) \leq 2 |\text{discr}(\Lambda)|$. 

\hfill

\begin{proof}
Let $w_0$ be a vector with minimal positive intersection $q(w_0,v)$. 
Then by Lemma 3.7. in \cite{_KV:fibrations_}, $q(w_0, v)$ divides 
$N = |\text{discr}(\Lambda)|$ (indeed, since $v$ is primitive, we can 
complete $v_1=v$ to a basis $\{v_1, \dots , v_r\}$ of $\Lambda$, 
and then $q(w_0, v) \Z$ is an ideal generated by $\{q(v, v_i)\}$, 
therefore the first column of the matrix $(q(v_j,v_i))$ is divisible 
by $q(w_0, v)$). Therefore, $0 < q(w_0, v) \leq N$. 
Let $\alpha$ be the smallest integer such that 
$q(w_0 + \alpha v, w_0 + \alpha v) >0$. Since $q(v,v)=0$, we have 
$q(w_0 + \alpha v, w_0 + \alpha v) = q(w_0,w_0) + 2 \alpha q(w_0, v)$. 
Then $w=w_0+\alpha v$ is a positive vector with $0 < q(w,v) = q(w_0,v) \leq N$. 
Notice that automatically $0<q(w,w) = q(w_0 + \alpha v, w_0 + \alpha v) = 
q(w_0,w_0) + 2 \alpha q(w_0, v) \leq 2N = 2 |\text{discr}(\Lambda)|$. 
\end{proof}

\hfill






One of our main results in \cite{_Kamenova_} was the following finiteness 
theorem with the assumption of existence of a Lagrangian fibration. Here 
we have dropped that assumption and we obtain a more general result. 

\hfill

\theorem\label{_finiteness_Theorem_1_}
There are at most finitely many deformation classes of hyperk\"ahler 
manifolds $M$ of dimension $2n$ with a fixed Fujiki constant $c$ and a given 
discriminant of the Beauville-Bogomolov-Fujiki lattice $(H^2(X,\Z),q)$. 

\hfill

\begin{proof}
If $b_2(M) = 3$ or $4$, there are finitely many lattices of this rank with 
fixed discriminant (see Theorem 1.1, chapter 9 of \cite{_Cassels_}). 
For each one of these finitely many lattices the positive vector $w$ with 
minimal positive intersection $q(w,w)$ has bounded self-intersection and we 
apply Charles's \ref{charles} directly. 

Now assume that $b_2 \geq 5$. By Meyer's theorem on indefinite lattices 
of rank at least $5$, there exists a non-trivial isotropic vector $v$ 
(\cite{_Cassels_}). We apply \ref{_bounded_lemma_} for 
$(\Lambda,q) = (H^2(X,\Z),q)$ and $v=v$. Then there exists a positive 
vector $w$ such that $0 < q(w,w) \leq 2|\text{discr}(\Lambda)|$. 
By Fujiki's formula, 
$$0< w^{2n} = c \cdot q(w,w)^n \leq c \cdot (2|\text{discr}(\Lambda)|)^n,$$
i.e., $w^{2n}$ is bounded by the given invariants. 
Deform $M$ to a hyperk\"ahler manifold $M'$ with a line bundle $L$ 
corresponding to the class $w$. For each top intersection 
$c_1(L)^{2n}$ in the interval $(0, (2|\text{discr}(\Lambda)|)^n]$, 
we obtain only finitely many deformation classes of $M$ by 
Charles's \ref{charles}. 
\end{proof}

\hfill

Since the families of hyperk\"ahler manifolds above 
form a bounded family, there are only finitely many 
choices of the second Betti number. The lattice $(H^2(M,\Z),q)$ encodes 
substantial information for hyperk\"ahler manifolds and $b_2(M)$ is an 
important invariant, therefore we list separately the following direct 
corollary. 

\hfill

\corollary
In the assumptions of \ref{_finiteness_Theorem_1_}, i.e., given numbers 
$n \in \Z_+$, $c \in \Q$ and discriminant $d \in \Z$, the second Betti 
number $b_2(M)$ is bounded. 

\hfill

We would like to mention also the following theorem in the recent literature. 

\hfill

\begin{theorem}(Sawon, \cite{_Sawon_finit_})
\label{_Sawon_finiteness_}
Fix positive integers $n$ and $d_1,\cdots,d_n$, with $d_1|d_2|\cdots |d_n$. 
Consider Lagrangian fibrations $\pi:M\rightarrow\C P^n$ that satisfy: 

(1) $\pi:M\rightarrow\C P^n$ admits a global section;

(2) there is a very ample line bundle on $M$ which gives a polarization of 
type $(d_1,\cdots,d_n)$ when restricted to a generic smooth fiber $M_t$;

(3) over a generic point $t$ of the discriminant locus the fiber 
$M_t$ is a rank-one semi-stable degeneration of abelian varieties;

(4) a neighbourhood $U$ of a generic point $t\in\C P^n$ describes a maximal 
variation of abelian varieties.

Then there are finitely many such Lagrangian fibrations up to deformation.
\end{theorem}

\hfill

Justin Sawon's proof is based on the following observations. 
The existence of a section implies that there is a distinguished point 
in each fiber. Together with the given type of a polarization it 
gives a natural classifying map 
$\phi : \C P^n \setminus \Delta \arrow {\cal A}_{d_1,\cdots,d_n}$, where 
$\Delta$ is the discriminant locus of $\pi$. Sawon extends $\phi$ to 
$\bar \phi : \C P^n \setminus \Delta_0 \arrow {\cal A}^*_{d_1,\cdots,d_n}$,
where ${\cal A}^*_{d_1,\cdots,d_n}$ is a partial compactification of 
${\cal A}_{d_1,\cdots,d_n}$ and $\Delta_0 \subset \Delta$ is of codimension 
$2$ in $\C P^n$. He chooses an ample line bundle $H$ on 
${\cal A}^*_{d_1,\cdots,d_n}$ and bounds $\deg (\bar \phi^* H)$, which is 
non-zero by $(4)$. Thus, $\bar \phi$ belongs to finitely many 
families of rational maps $\C P^n \dashrightarrow {\cal A}^*_{d_1,\cdots,d_n}$, 
and finiteness of Lagrangian fibrations as above follows. 

\hfill

\remark
Notice that as a corollary of Matsushita's conjecture, part $(4)$ of 
Sawon's Theorem simply excludes locally isotrivial fibrations. We need to 
apply only the deformational version (see van Geemen-Voisin's \ref{gv-m}) 
of Matsushita's conjecture to Sawon's 
theorem in order to remove the seemingly restrictive assumption $(4)$. 

\hfill 

\remark
If there is a section $\sigma : \C P^n \arrow M$, this means 
that $\sigma (\C P^n)$ is a Lagrangian subvariety in $M$. 
Finding Lagrangian subvarieties in a hyperk\"ahler manifold is itself a 
very interesting problem (for example, see \cite{_H_T_}).

\hfill

Using the methods above combined with F. Charles' \ref{charles}, 
in \cite{_Kamenova_}  we generalized Sawon's \ref{_Sawon_finiteness_}. 

\hfill

\theorem \label{_finiteness_Theorem_2_}
Consider a Lagrangian fibration $\pi:M \arrow\C P^n$ such that 
there is a line bundle $P$ on $M$ with $q(P)>0$ and with a given $P$-degree 
$d$ on the general fiber $F$ of $\pi$, i.e., $P^n \cdot F = d$. 
Then there are at most finitely many deformation classes of hyperk\"ahler 
manifolds $M$ as above, i.e., they form a bounded family. 

\hfill

\begin{proof}
Let $L = \pi^* {\cal O} (1)$ be a nef parabolic class ($q(L)=0$) coming 
from the Lagrangian fibration. 
The fundamental class of the general fiber $F$ of $\pi$ is $[F] = L^n$. 
By assumption, $P^n \cdot L^n = d$ is fixed. 
Define $v = L/m$, where $m \in \Z_{>0}$ is the divisibility of $L$, 
and therefore $v$ is a primitive class. 
Since $P$ is in the interior of the of the positive cone 
${\cal C}$ 
and $v$ is on the boundary of $\cal C$, it follows that $q(P, v)>0$ 
(Corollary 7.2 in \cite{_BHPV_}). 
Now we shall follow the proof of \ref{_bounded_lemma_}. 
Let $k$ be the smallest integer such that $q(P+kv)>0$. 
Then $q(P+kv) \leq 2 q(P,v)$ and by applying 
Fujiki's formula twice (as in \ref{Fujiki_formula} and \ref{_P_L_identity_}), 
we obtain: 
$$(P+kv)^{2n} = c \cdot q(P+kv)^n \leq c 2^n q(P,v)^n 
= {2n \choose n} P^n \cdot v^n =$$ $${2n \choose n} 
\frac{P^n \cdot L^n}{m^n} = {2n \choose n} \frac{d}{m^n} \leq 
{2n \choose n} d.$$
By Charles's \ref{charles} there is a bounded family of such $M$, 
which implies finiteness of deformations of $M$.
\end{proof}

\hfill

\remark
In \ref{_finiteness_Theorem_2_} we proved 
finiteness of deformation classes of the total space $M$ of the Lagrangian 
fibration. However, in \ref{KVfin} the author together with Misha Verbitsky 
proved that for a fixed compact manifold $M$ there are only finitely many 
deformation types of hyperk\"ahler Lagrangian fibrations with total space $M$ 
provided that $b_2(M) \geq 7$. For all known examples of hyperk\"ahler 
manifolds one has $b_2(M) \geq 7$ and it is suspected that this is 
always the case. 

\hfill

\noindent{\bf Acknowledgments.} 
The author is very grateful to the organizers of the Simons Semester in Poland 
for their kind invitation and to the Banach Center in B\c{e}dlewo for their 
excellent hospitality and conference. She would like to thank 
Fran\c{c}ois Charles, Misha Verbitsky, Dennis Sullivan, Ulrike Rie{\ss} and 
Michel Brion for interesting conversations and suggestions.


{\small

\hfill

\noindent {\sc Ljudmila Kamenova\\
Department of Mathematics, 3-115 \\
Stony Brook University \\
Stony Brook, NY 11794-3651, USA,} \\
\tt kamenova@math.stonybrook.edu\\

\end{document}